\documentclass{llncs}

\usepackage{t1enc}
\usepackage[latin2]{inputenc}
\usepackage[english]{babel}
\usepackage{amsmath}
\usepackage{amssymb}
\usepackage{bm}
\usepackage[all]{xy}
\usepackage[mathscr]{eucal}
\usepackage{xspace}
\newcommand{\G}{\mathscr{G}}

\newcommand{\RICA}{\mbox{$\mathbb{R}$-ICA}\xspace}
\newcommand{\RISA}{\mbox{$\mathbb{R}$-ISA}\xspace}
\newcommand{\CICA}{\mbox{$\mathbb{C}$-ICA}\xspace}
\newcommand{\CISA}{\mbox{$\mathbb{C}$-ISA}\xspace}
\newcommand{\KICA}{\mbox{$\mathbb{K}$-ICA}\xspace}
\newcommand{\KISA}{\mbox{$\mathbb{K}$-ISA}\xspace}

\newcommand{\K}{\ensuremath{\mathbb{K}}\xspace}
\newcommand{\R}{\ensuremath{\mathbb{R}}\xspace}
\newcommand{\C}{\ensuremath{\mathbb{C}}\xspace}

\newcommand{\pv}{\ensuremath{\varphi_{v}}}
\newcommand{\pM}{\ensuremath{\varphi_{M}}}

\title{Separation Theorem for \K-Independent Subspace Analysis with Sufficient Conditions}
\author{Zolt\'an Szab\'o, Barnab\'as P\'oczos, and Andr\'as L{\H o}rincz}
\institute{Department of Information Systems, E\"{o}tv\"{o}s Lor{\'a}nd University,\\
Pázmány P. sétány 1/C, Budapest H-1117, Hungary\\
Research Group on Intelligent Information Systems \\ Hungarian Academy of
Sciences\\
\email{szzoli@cs.elte.hu, pbarn@cs.elte.hu, lorincz@inf.elte.hu},\\
WWW home page: \url{http://nipg.inf.elte.hu}}

\begin{document}
\maketitle
\begin{abstract}
Here, a Separation Theorem about \K-Independent Subspace Analysis ($\K\in\{\R,\C\}$ real or complex), a generalization
of \K-Independent Component Analysis (\KICA) is proven. According to the theorem, \KISA estimation can be executed in
two steps under certain conditions. In the first step, 1-dimensional \KICA estimation is executed. In the second step,
optimal permutation of the \KICA elements is searched for. We present sufficient conditions for the \KISA Separation
Theorem. Namely, we shall show that (i) spherically symmetric sources (both for real and complex cases), as well as
(ii) real 2-dimensional sources invariant to 90$^\circ$ rotation, among others, satisfy the conditions of the theorem.
\end{abstract}

\section{Introduction}
(Real) Independent Component Analysis (\RICA) \cite{jutten91blind,comon94independent} aims to recover linearly or non-linearly mixed
independent and hidden sources. There is a broad range of applications for \RICA, such as blind source separation and blind source
deconvolution \cite{bell95information}, feature extraction \cite{bell97independent}, denoising \cite{hyvarinen99sparse}. Particular
applications include, e.g., the analysis of financial data \cite{kiviluoto98independent}, data from neurobiology, fMRI, EEG, and MEG (see,
e.g., \cite{makeig96independent,Vigario98independent} and references therein). For a recent review on \RICA see \cite{choi05blind}.

Original \RICA algorithms are 1-dimensional in the sense that all sources are assumed to be independent real valued random variables.
However, applications where not all, but only certain groups of the sources are independent may have high relevance in practice. In this
case, independent sources can be multi-dimensional. For example, consider the generalization of the cocktail-party problem, where
independent groups of people are talking about independent topics, or that more than one group of musicians are playing at the party. The
separation task requires an extension of \RICA, which can be called (Real) Independent Subspace Analysis (\RISA)
\cite{hyvarinen00emergence}, Multi-Dimensional Independent Component Analysis (MICA) \cite{cardoso98multidimensional}, Group ICA
\cite{theis05blind}, and Independent Vector Analysis (IVA) \cite{kim06independent}. Throughout the paper, we shall use the first
abbreviation. An important application for \RISA is, e.g., the processing of EEG-fMRI data \cite{akaho99MICA}.

Efforts have been made to develop \RISA algorithms
\cite{cardoso98multidimensional,akaho99MICA,vollgraf01multi,bach03finding,poczos05independent1,poczos05independent2,theis05blind}.
Related theoretical problems concern mostly the estimation of entropy or mutual information. In this context, entropy
estimation by Edgeworth expansion \cite{akaho99MICA} has been extended to more than 2 dimensions and has been used for
clustering and mutual information testing \cite{hulle05edgeworth}. $k$-nearest neighbors and geodesic spanning tree
methods have been applied in \cite{poczos05independent1} and \cite{poczos05independent2} for the \RISA problem. Other
recent approaches search for independent subspaces via kernel methods \cite{bach03finding} and joint block
diagonalization \cite{theis05blind}.

Beyond the case of real numbers, the search for complex components (Complex Independent Component Analysis, \CICA)
assumes more and more practical relevance. Such problems include, beyond others, (i) communication systems, (ii)
biomedical signal processing, e.g., processing of (f)MRI and EEG data, brain modelling, (iii) radar applications, (iv)
frequency domain methods (e.g., convolutive models).

There is a large number of existing \CICA procedures
\cite{anemuller03complex,back94blind,fiori00blind,bingham00fast,comon94independent,cardoso93blind,cardoso99highorder,moreau01generalization,cardoso96equivariant,fiori02complex,belouchrani97blind,arie00blind,smaragdis98blind,calhoun02complex,xu04minimax,eriksson04complex,douglas06equivariant,douglas06fixed}.
Maximum likelihood principle and complex recurrent neural network are used in \cite{anemuller03complex}, and
\cite{back94blind}, respectively. The APEX algorithm is based on Hebbian learning \cite{fiori00blind}. Complex FastICA
algorithm can be found in \cite{bingham00fast}. More solutions are based on cumulants: e.g., \cite{comon94independent},
the JADE algorithm \cite{cardoso93blind,cardoso99highorder}, its higher order variants \cite{moreau01generalization},
and the EASI algorithm family \cite{cardoso96equivariant}. `Rigid-body' learning theory is used in
\cite{fiori02complex}. The SOBI algorithm \cite{belouchrani97blind} searches for joint diagonalizer matrix, its refined
version, the WASOBI method \cite{arie00blind} approximates by means of weighted nonlinear least squares. There are
complex variants of the infomax technique, such as the a split-complex \cite{smaragdis98blind} and the fully-complex
infomax \cite{calhoun02complex} procedures. Minimax Mutual Information \cite{xu04minimax} and strong-uncorrelating
transforms \cite{eriksson04complex,douglas06equivariant,douglas06fixed} make further promising directions.

An important observation of previous computer studies \cite{cardoso98multidimensional,poczos05independent3} is that
general \RISA solver algorithms are not more efficient, in fact, sometimes produce lower quality results than simple
\RICA algorithm superimposed with searches for the optimal permutation of the components. This observation led to the
present theoretical work and to some computer studies that have been published elsewhere
\cite{szabo06cross,szabo06real}. We treat both the real and the complex cases.

This technical report is constructed as follows: Section~\ref{sec:preliminaries} introduces complex random variables. In
Section~\ref{sec:ISA-model} the \KISA task is described. Section~\ref{sec:KISA-sep-theorem} contains our Separation Theorem for the \KISA
task. Sufficient conditions for the theorem are provided in Section~\ref{sec:suff-cond}. Conclusions are drawn in
Section~\ref{sec:conclusion}.

\section{Basic Concepts: Matrices, Complex Random Variables}\label{sec:preliminaries}
We introduce the basic concepts for using complex random variables. Excellent review can be found in \cite{eriksson06complex}.

$\mathbf{B}^T$ is the transposed of matrix $\mathbf{B}\in \C^{L\times L}$. Complex conjugation is denoted by a bar it
concerns all elements of a matrix. The transposed complex conjugate of matrix $\mathbf{B}$ is the adjoint matrix
$\mathbf{B}^*=\bar{\mathbf{B}}^T$. Matrix $\mathbf{B}\in\C^{L\times L}$ is called \emph{unitary} if
$\mathbf{B}\mathbf{B}^*=\mathbf{I}_L$, \emph{orthogonal} if $\mathbf{B}\mathbf{B}^T=\mathbf{I}_L$, where $\mathbf{I}_L$
is the \mbox{$L$-dimensional} identity matrix. The sets of $L\times L$ dimensional unitary and orthogonal matrices are
denoted by $\mathcal{U}^D$ and $\mathcal{O}^L$, respectively.

A \emph{complex-valued random variable} $\mathbf{u}\in\C^L$ (shortly complex random variable) is defined as a random
variable of the form \mbox{$\mathbf{u}=\mathbf{u}_R+i\mathbf{u}_I$}, where the real and imaginary parts of
$\mathbf{u}$, i.e., $\mathbf{u}_R$ and $\mathbf{u}_I\in\R^L$ are real random variables, $i=\sqrt{-1}$. Expectation
value of complex random variables is \mbox{$E[\mathbf{u}]=E[\mathbf{u}_R]+iE[\mathbf{u}_I]$}, and the variable can be
characterized in second order by its \emph{covariance matrix}
$cov[\mathbf{u}]=E[(\mathbf{u}-E[\mathbf{u}])(\mathbf{u}-E[\mathbf{u}])^*]$ and by its \emph{pseudo-covariance matrix}
\mbox{$pcov[\mathbf{u}]=E[(\mathbf{u}-E[\mathbf{u}])(\mathbf{u}-E[\mathbf{u}])^T]$}. Complex random variable
$\mathbf{u}$ is called \emph{full}, if $cov[u]$ is positive definite. Throughout this paper all complex variables are
assumed to be full (that is, they are not concentrated in any lower dimensional complex subspace).

\section{The \KISA Model}\label{sec:ISA-model}
First, Section~\ref{subsec:KISA-eqs} introduces the \KISA task. Section~\ref{subsec:CISA ambiguities} is about the
ambiguities of the problem. Section~\ref{subsec:KISA-cost-function} defines the entropy based cost function of the
\KISA task.

\subsection{The \KISA Equations}\label{subsec:KISA-eqs}
We define the complex (and real) ISA task. The two models are treated together by using the joint notation
$\K\in\{\R,\C\}$. Assume that we have $M$ hidden and independent \emph{components} (random variables) and that only the
mixture of them is available for observation:
\begin{eqnarray}
\mathbf{z}(t)&=&\mathbf{A}\mathbf{s}(t),\label{eq:M2}
\end{eqnarray}
where $\mathbf{s}(t)=\left[\mathbf{s}^1(t);\ldots;\mathbf{s}^M(t)\right]$ is a vector concatenated of components
$\mathbf{s}^m(t)\in\K^{d}$. For fix $m$, $\mathbf{s}^m(t)$ is i.i.d. (independent identically distributed) in $t$,
$\mathbf{s}^i$ is independent from $\mathbf{s}^j$, if $i\neq j$. The total dimension of the components is
\mbox{$D:=Md$}. Thus, \mbox{$\mathbf{s}(t), \mathbf{z}(t)\in\K^D$}. Matrix \mbox{$\mathbf{A}\in\K^{D\times D}$} is the
so called \emph{mixing matrix} which, according to our assumptions, is invertible.

The goal of the \KISA problem is to estimate  the original source $\mathbf{s}(t)$ and the unknown mixing matrix
$\mathbf{A}$ (or its inverse $\mathbf{W}$, which is called the \emph{separation matrix}) by using observations
$\mathbf{z}(t)$ only. We talk about complex ISA task if $\K=\C$ (\CISA), and real ISA task if $\K=\R$ (\RISA). If
$d=1$, then complex ICA (\CICA), and real ICA (\RICA) tasks are obtained.

\subsection{Ambiguities of the \KISA Model}\label{subsec:CISA ambiguities}
Identification of the $\KISA$ model is ambiguous. However, there are obvious ambiguities of the model: hidden
components can be determined up to permutation of subspaces and invertible transformation within the subspaces. Further
details concerning ambiguities can be found here: \RICA \cite{eriksson04contributions}, \RISA
\cite{theis04uniqueness1}, \CICA
\cite{eriksson04contributions,eriksson06complex,theis04uniqueness1,theis04uniqueness2}, \CISA (see
Appendix~\ref{sec:ISA separability}).

Ambiguities within subspaces can be lessened. Namely, given our assumption on the invertibility of matrix $\mathbf{A}$, we can assume
without any loss of generality that both the sources and the observation are \emph{white}, that is,
\begin{eqnarray}
E[\mathbf{s}]&=&\mathbf{0},cov\left[\mathbf{s}\right]=\mathbf{I}_D,\label{eq:white1}\\
E[\mathbf{z}]&=&\mathbf{0},cov\left[\mathbf{z}\right]=\mathbf{I}_D\label{eq:white2}.
\end{eqnarray}
Below, we treat real and complex cases separately:

\paragraph{Real case:} It then follows that the mixing matrix $\mathbf{A}$ and thus the separation matrix $\mathbf{W}=\mathbf{A}^{-1}$ are
orthogonal:
\begin{equation}
\mathbf{I}_D
=cov[\mathbf{z}]=E\left[\mathbf{z}\mathbf{z}^T\right]=\mathbf{A}E\left[\mathbf{s}\mathbf{s}^T\right]\mathbf{A}^T=\mathbf{A}\mathbf{I}_D\mathbf{A}^T=\mathbf{A}\mathbf{A}^T.
\end{equation}
The ambiguity of the ISA task is decreased by Eqs.~\eqref{eq:white1}--\eqref{eq:white2}: Now, $\mathbf{s}^m$ sources
are determined up to permutation \emph{and} orthogonal transformation.

\paragraph{Complex case:}
It then follows that the mixing matrix $\mathbf{A}$ and thus the separation matrix $\mathbf{W}=\mathbf{A}^{-1}$ are
unitary:
        \begin{equation}
        \mathbf{I}_D =cov[\mathbf{z}]=E\left[\mathbf{z}\mathbf{z}^*\right]=\mathbf{A}E\left[\mathbf{s}\mathbf{s}^*\right]\mathbf{A}^*=\mathbf{A}\mathbf{I}_D\mathbf{A}^*=\mathbf{A}\mathbf{A}^*.
        \end{equation}
Thus, components $\mathbf{s}^m$ are determined up to permutation \emph{and} unitary transformation within the subspace.

\subsection{The \KISA Cost Function}\label{subsec:KISA-cost-function}
Now we sketch how to transcribe the \KISA task into the minimization of sum of multi-dimensional entropies for
orthogonal matrices (in the real case) and for unitary matrices (in the complex case). We shall use these formulations
of the \KISA task to prove the real and complex versions of the Separation Theorem
(Section~\ref{sec:KISA-sep-theorem}).

\subsubsection{Real Case}
The \RISA task can be viewed as the minimization of mutual information between the estimated components:
\begin{equation}
 I\left(\mathbf{y}^1,\ldots,\mathbf{y}^M\right):=\int
                f(\mathbf{v})\ln\left[\frac{f(\mathbf{v})}{\prod_{i=1}^Mf_m(v_m)}\right]\mathrm{d}\mathbf{v}\label{eq:mutualinfo}
\end{equation}
on the orthogonal group ($\mathbf{W}\in\mathcal{O}^D$), where \mbox{$\mathbf{y}=\mathbf{W}\mathbf{z}$},
\mbox{$\mathbf{y}=\left[\mathbf{y}^1;\ldots;\mathbf{y}^M\right]$}, $f$ and $f_m$ are density functions of $\mathbf{y}$
and marginals $\mathbf{y}^m$, respectively. This cost function $I$ is equivalent to the minimization of the sum of
d-dimensional entropies, because
\begin{eqnarray}
I\left(\mathbf{y}^1,\ldots,\mathbf{y}^M\right)&=&\sum_{m=1}^MH\left(\mathbf{y}^m\right)-H(\mathbf{y})\label{eq:I2H}\\
 &=&\sum_{m=1}^MH\left(\mathbf{y}^m\right)-H(\mathbf{Wz})\\
 &=&\sum_{m=1}^MH\left(\mathbf{y}^m\right)-[H(\mathbf{z})+\ln(\left|\det(\mathbf{W})\right|)].
\end{eqnarray}

Here, $H$ is Shannon's (multi-dimensional) differential entropy defined with logarithm of base $e$,
$\left|\cdot\right|$ denotes absolute value, `$\det$' stands for determinant. In the second equality, the
$\mathbf{y}=\mathbf{W}\mathbf{z}$ relation was exploited, and the
\begin{equation}
H(\mathbf{Wz})=H(\mathbf{z})+\ln\left(\left|\det(\mathbf{W})\right|\right)\label{eq:R-entropy-trafo}
\end{equation}
rule describing transformation of the differential entropy \cite{cover91elements} was used. $\det(\mathbf{W})=1$
because of the orthogonality of $\mathbf{W}$, so $\ln(\left|\det(\mathbf{W})\right|)=0$. The $H(\mathbf{z})$ term of
the cost is constant in $\mathbf{W}$, therefore the \RISA task is equivalent to the minimization of the cost function
\begin{equation}
J(\mathbf{W}):=\sum_{m=1}^MH\left(\mathbf{y}^m\right)\label{eq:ISA-costfunction}\rightarrow\min_{\mathbf{W}\in\mathcal{O}^D}.
\end{equation}
\subsubsection{Complex Case}
Similarly, the \CISA task can be viewed as the minimization of mutual information between the estimated components [see
Eq.~\eqref{eq:mutualinfo}], but on the unitary group ($\mathbf{W}\in\mathcal{U}^D$). Here, the Shannon entropy of
random variable $\C^L\ni\mathbf{u}$ $(\mathbf{y}^m$, or $\mathbf{y})$ is the entropy of $\pv(\mathbf{u})\in\R^{2L}$,
where
\begin{equation}
    \pv:\C^L\ni\mathbf{u}\mapsto\mathbf{u}\otimes\left[\begin{array}{c}\Re(\cdot)\\\Im(\cdot)\end{array}\right]\in\R^{2L}.\label{eq:pv}
\end{equation}
That is, $H(\mathbf{u}):=H[\pv(\mathbf{u})]$. Here: $\otimes$ is the Kronecker product, $\Re$ stands for the real part,
$\Im$ for the imaginary part, subscript 'v' for vector. One can neglect the last term of the $H(\mathbf{y})$  cost
function [see, Eq.~\eqref{eq:I2H}] during optimization (alike in the real case). To see this, consider the mapping
\begin{equation}
    \pM:\C^{L\times L}\ni\mathbf{M}\mapsto \mathbf{M}\otimes\left[\begin{array}{rr}\Re(\cdot)&-\Im(\cdot)\\\Im(\cdot)&\Re(\cdot)\end{array}\right]\in\R^{2L\times
    2L},\label{eq:pM}
\end{equation} where subscript `$M$' indicate matrices. Known properties of mappings $\pv$, $\pM$ are as follows \cite{krishnaiah86complex}:
\begin{align}
    \det[\pM(\mathbf{M})]&=|\det(\mathbf{M})|^2\label{prop:det-trafo},\\
    \pv(\mathbf{M}\mathbf{v})&=\pM(\mathbf{M})\pv(\mathbf{v}).\label{prop:mtx-vec-hom}
\end{align}
In words: \eqref{prop:det-trafo} describes transformation of determinant, \eqref{prop:mtx-vec-hom} expresses
preservation of operation for matrix-vector multiplication.\footnote{~Note that this connection allows one to reduce
the \CISA task (and thus the \CICA task) to a \RISA task directly. According to our experiences, however, methods that
rely on the \C-Separation Theorem that we present here are much more efficient.} The following relation holds for the
entropy transformation of complex variables:
\begin{lemma}[Transformation of entropy for complex variables]\label{lem:C-H-trafo}
Let $\mathbf{u}\in\C^L$ denote a random variable and let $\mathbf{V}\in\C^{L\times L}$ be a matrix. Then
\begin{equation}
     H(\mathbf{Vu})=H(\mathbf{u})+\ln\left(\left|\det(\mathbf{V})\right|^2\right)\label{eq:H(W*C vv)}
\end{equation}
\end{lemma}
\begin{proof}
\begin{align}
H(\mathbf{V}\mathbf{u})&=H[\pv(\mathbf{V}\mathbf{u})]=H[\pM(\mathbf{V})\pv(\mathbf{u})]=H[\pv(\mathbf{u})]+\ln(|\det[\pM(\mathbf{V})]|)
=
H[\pv(\mathbf{u})]+\ln(|\det\mathbf{V}|^2)\\
&= H(\mathbf{u})+\ln(|\det\mathbf{V}|^2)
\end{align}
The above steps can be justified as follows:
\begin{enumerate}
    \item
        the first equation uses the definition of entropy for complex variables,
    \item
       then we used property \eqref{prop:mtx-vec-hom},
    \item
        transformed the entropy of random variables in $\R^{2L}$ [see, Eq.~\eqref{eq:R-entropy-trafo}].
    \item
        exploited \eqref{prop:det-trafo}, and
    \item
       applied the definition of entropy for complex variables again.
\end{enumerate}\qed
\end{proof}

Thus $H(\mathbf{y})=H(\mathbf{Wz})+\ln(1)=H(\mathbf{z})$, where unitarity of $\mathbf{W}$ is exploited. Further,
$H(\mathbf{z})$ is not dependent of matrix $\mathbf{W}$ and thus term $H(\mathbf{y})$ can be neglected during the
course of optimization. We conclude that the $\C$-ISA task can be written as the minimization of sum of
multi-dimensional entropies. The cost function to be optimized within unitary matrices:
\begin{equation}
       J(\mathbf{W})=\sum_{m=1}^M H(\mathbf{y}^m)\rightarrow\min_{\mathbf{W}\in\mathcal{U}^D}.
\end{equation}

\section{The \RISA and \CISA Separation Theorem}\label{sec:KISA-sep-theorem}
The main result of this work is that the \KISA task may be accomplished in two steps under certain conditions. In the first step \KICA is
executed. The second step is search for the optimal permutation of the \KICA components. Section~\ref{subsec:RISA-sep-theorem} is about the
real, whereas Section~\ref{subsec:CISA-sep-theorem} is about the complex case.

\subsection{The \RISA Separation Theorem}\label{subsec:RISA-sep-theorem}
We shall rely on entropy inequalities (Section~\ref{subsubsec:REPIineq-s}). Connection to the \RICA cost function is derived in
Section~\ref{subsubsec:connection2RICA}. Finally, Section~\ref{subsubsec:RISA-sepT-proof} contains the proof of our theorem.

\subsubsection{EPI-type Relations (Real Case)}\label{subsubsec:REPIineq-s}
First, consider the so called Entropy Power Inequality (EPI)
\begin{equation}
e^{2H\left(\sum_{i=1}^Lu_i\right)}\ge \sum_{i=1}^L e^{2H(u_i)},\label{eq:EPI}
\end{equation}
where $u_1,\ldots,u_L\,\in\R$ denote continuous random variables (The name of this inequality is \R-EPI, because we shall need its complex
variant later). This inequality holds for example, for independent continuous variables \cite{cover91elements}.

Let $\left\|\cdot\right\|$ denote the Euclidean norm. That is, for $\mathbf{w}\in\R^L$
\begin{equation}
\left\|\mathbf{w}\right\|^2:=\sum_{i=1}^Lw_i^2,
\end{equation}
where $w_i$ is the $i^{th}$ coordinate of vector $\mathbf{w}$. The surface of the \mbox{$L$-dimensional} unit sphere
shall be denoted by $S^L(\R)$:
\begin{equation}
S^L(\R):=\{\mathbf{w}\in\R^L:\left\|\mathbf{w}\right\|=1\}.
\end{equation}
If \R-EPI is satisfied [on $S^L(\R)$] then a further inequality holds:
\begin{lemma}\label{lem:suff}
Suppose that continuous random variables $u_1,\ldots,u_L \,\in\R$ satisfy the following inequality
\begin{equation}
    e^{2H\left(\sum_{i=1}^Lw_iu_i\right)}\ge \sum_{i=1}^L e^{2H(w_iu_i)}, \forall\mathbf{w}\in S^L(\R).\label{eq:w-EPI}
\end{equation}
This inequality will be called the \emph{\R-w-EPI} condition. Then
\begin{equation}
H\left(\sum_{i=1}^L w_iu_i\right)\ge\sum_{i=1}^Lw_i^2H\left(u_i\right), \forall\mathbf{w}\in S^L(\R).\label{eq:suff}
\end{equation}
\end{lemma}

\begin{note}
\R-w-EPI holds, for example, for independent variables $u_i$, because independence is not affected by multiplication with a constant.
\end{note}

\begin{proof} Assume that $\mathbf{w}\in S^L(\R)$. Applying
$\ln$ on condition \eqref{eq:w-EPI}, and using the monotonicity of the $\ln$ function, we can see that the first inequality is valid in the
following inequality chain
\begin{equation}
2H\left(\sum_{i=1}^Lw_iu_i\right)\ge \ln\left(\sum_{i=1}^L e^{2H(w_iu_i)}\right)=\ln\left(\sum_{i=1}^Le^{2H(u_i)}\cdot
w_i^2\right)\ge\sum_{i=1}^Lw_i^2\cdot\ln\left(e^{2H(u_i)}\right)=\sum_{i=1}^Lw_i^2\cdot2H(u_i).
\end{equation}
Then,
\begin{enumerate}
    \item
        we used the relation \cite{cover91elements}:
        \begin{equation}
            H(w_iu_i)=H(u_i)+\ln\left(\left|w_i\right|\right)
        \end{equation}
        for the entropy of the transformed variable. Hence
        \begin{equation}
            e^{2H(w_iu_i)}=e^{2H(u_i)+2\ln\left(\left|w_i\right|\right)}=e^{2H(u_i)}\cdot
            e^{2\ln\left(\left|w_i\right|\right)}=e^{2H(u_i)}\cdot
            w_i^2.\label{eq:2entr-transf}
        \end{equation}
    \item
        In the second inequality, we exploited the concavity of $\ln$.\qed
\end{enumerate}
\end{proof}

\subsubsection{Connection to the Cost Function of the \RICA Task}\label{subsubsec:connection2RICA} Now we shall use Lemma~\ref{lem:suff} to proceed. The
\RISA Separation Theorem will be a corollary of the following claim:
\begin{proposition}\label{R-prop}
Let $\mathbf{y}=\left[\mathbf{y}^1;\ldots;\mathbf{y}^M\right]=\mathbf{y}(\mathbf{W})=\mathbf{W}\mathbf{s}$, where
$\mathbf{W}\in \mathcal{O}^D$, $\mathbf{y}^m$ is the estimation of the $m^{th}$ component of the \RISA task. Let
$y^m_i$ be the $i^{th}$ coordinate of the $m^{th}$ component. Similarly, let $s^m_i$ stand for the $i^{th}$ coordinate
of the $m^{th}$ source. Let us assume that the $\mathbf{s}^m$ sources satisfy condition~\eqref{eq:suff}. Then
\begin{equation}\label{eq:main-prop}
\sum_{m=1}^M\sum_{i=1}^dH\left(y^m_i\right)\ge
\sum_{m=1}^M\sum_{i=1}^dH\left(s^m_i\right).
\end{equation}
\end{proposition}

\begin{proof}
Let us denote the $(i,j)^{th}$ element of matrix $\mathbf{W}$ by $W_{i,j}$. Coordinates of $\mathbf{y}$ and
$\mathbf{s}$ will be denoted by $y_i$ and $s_i$, respectively. Further, let $\G^1, \ldots, \G^M$ denote the indices of
the $1^{st}, \ldots , M^{th}$ subspaces, i.e., $\G^1:=\{1,\ldots,d\},\ldots,\G^M:=\{D-d+1,\ldots,D\}$. Now, writing the
elements of the $i^{th}$ row of matrix multiplication $\mathbf{y}=\mathbf{W}\mathbf{s}$, we have
\begin{equation}
y_i=\sum_{j\in \G^1} W_{i,j}s_j+\ldots+\sum_{j\in \G^M} W_{i,j}s_j\label{eq:y=Ws}
\end{equation}
and thus,
\begin{eqnarray}
\lefteqn{H\left(y_i\right)=}\nonumber\\
&&=H\left(\sum_{j\in \G^1} W_{i,j}s_j+\ldots+\sum_{j\in \G^M} W_{i,j}s_j\right)\label{eq:H(y=Ws)}\\
&&=H\left(\left(\sum_{l\in\G^1}W_{i,l}^2\right)^{\frac{1}{2}}\frac{\sum_{j\in\G^1}W_{i,j}s_j}{\left(\sum_{l\in\G^1}W_{i,l}^2\right)^{\frac{1}{2}}}
+ \ldots + \left(\sum_{l\in\G^M}W_{i,l}^2\right)^{\frac{1}{2}}\frac{\sum_{j\in\G^M}W_{i,j}s_j}{\left(\sum_{l\in\G^M}W_{i,l}^2\right)^{\frac{1}{2}}}\right)\label{eq:w^2-in}\\
&&\ge\left(\sum_{l\in\G^1}W_{i,l}^2\right)H\left(\frac{\sum_{j\in\G^1}W_{i,j}s_j}{\left(\sum_{l\in\G^1}W_{i,l}^2\right)^{\frac{1}{2}}}\right)
+ \ldots    + \left(\sum_{l\in\G^M}W_{i,l}^2\right) H\left(\frac{\sum_{j\in\G^M}W_{i,j}s_j}{\left(\sum_{l\in\G^M}W_{i,l}^2\right)^{\frac{1}{2}}}\right)\label{eq:Lem2-applied}\\
&&=\left(\sum_{l\in\G^1}W_{i,l}^2\right)
H\left(\sum_{j\in\G^1}\frac{W_{i,j}}{\left(\sum_{l\in\G^1}W_{i,l}^2\right)^{\frac{1}{2}}}s_j\right)
+ \ldots + \left(\sum_{l\in\G^M}W_{i,l}^2\right) H\left(\sum_{j\in\G^M}\frac{W_{i,j}}{\left(\sum_{l\in\G^M}W_{i,l}^2\right)^{\frac{1}{2}}}s_j\right)\label{eq:Lem2-applied-again-pre}\\
&&\ge\left(\sum_{l\in\G^1}W_{i,l}^2\right)
\sum_{j\in\G^1}\left(\frac{W_{i,j}}{\left(\sum_{l\in\G^1}W_{i,l}^2\right)^{\frac{1}{2}}}\right)^2H\left(s_j\right)\label{eq:Lem2-applied-again}
+ \ldots + \left(\sum_{l\in\G^M}W_{i,l}^2\right) \sum_{j\in\G^M}\left(\frac{W_{i,j}}{\left(\sum_{l\in\G^M}W_{i,l}^2\right)^{\frac{1}{2}}}\right)^2H\left(s_j\right)\\
&&=\sum_{j\in\G^1}W_{i,j}^2H\left(s_j\right)+\ldots+\sum_{j\in\G^M}W_{i,j}^2H\left(s_j\right)\label{eq:H(yi-last)}
\end{eqnarray}
The above steps can be justified as follows:
\begin{enumerate}
    \item
        \eqref{eq:H(y=Ws)}: Eq.~\eqref{eq:y=Ws} was inserted into the argument of $H$.
    \item
        \eqref{eq:w^2-in}: New terms were added for Lemma~\ref{lem:suff}.
    \item
        \eqref{eq:Lem2-applied}: Sources $\mathbf{s}^m$ are independent of each other and this independence is preserved upon
        mixing \emph{within} the subspaces, and we could also use Lemma~\ref{lem:suff}, because $\mathbf{W}$
        is an orthogonal matrix.
    \item
        \eqref{eq:Lem2-applied-again-pre}: Nominators were transferred into the $\sum_j$ terms.
    \item
        \eqref{eq:Lem2-applied-again}: Variables $\mathbf{s}^m$ satisfy condition~\eqref{eq:suff} according
        to our assumptions.
    \item
        \eqref{eq:H(yi-last)}: We simplified the expression after squaring.
\end{enumerate}
Using this inequality, summing it for $i$, exchanging the order of the sums, and making use of the orthogonality of
matrix $\mathbf{W}$, we have
\begin{eqnarray}
\sum_{i=1}^DH(y_i)&\ge&\sum_{i=1}^D\left(\sum_{j\in\G^1}W_{i,j}^2H\left(s_j\right)+\ldots+\sum_{j\in\G^M}W_{i,j}^2H\left(s_j\right)\right)\\
&=&\sum_{j\in\G^1}\left(\sum_{i=1}^DW^2_{i,j}\right)H\left(s_j\right)+\ldots+\sum_{j\in\G^M}\left(\sum_{i=1}^DW^2_{i,j}\right)H\left(s_j\right)\\
&=&\sum_{j=1}^DH(s_j).
\end{eqnarray}
\qed
\end{proof}

\begin{note}
The proof holds for components with different dimensions. This is also true for the following theorem.
\end{note}

\subsubsection{Proof of the $\R$-ISA Separation Theorem}\label{subsubsec:RISA-sepT-proof}

Having this proposition, now we present our main theorem.
\begin{theorem}[Separation Theorem for \RISA]
Presume that the $\mathbf{s}^m$ sources of the \RISA model satisfy condition~\eqref{eq:suff}, and that the \RICA cost
function $J(\mathbf{W})=\sum_{m=1}^M\sum_{i=1}^dH(y^m_i)$ has minimum ($\mathbf{W}\in\mathscr{O}^D$). Then it is
sufficient to search for the minimum of the \RISA task as a permutation of the solution of the \RICA task. Using the
concept of separation matrices, it is sufficient to explore forms
\begin{equation}
    \mathbf{W}_{\text{\RISA}}=\mathbf{P}\mathbf{W}_{\text{\RICA}},
\end{equation}
where $\mathbf{P} \left(\in\R^{D\times D}\right)$ is a permutation matrix to be determined.
\end{theorem}
\begin{proof}
\RICA minimizes the l.h.s. of Eq.~\eqref{eq:main-prop}, that is, it minimizes \mbox{$\sum_{m=1}^M\sum_{i=1}^dH\left(y^m_i\right)$}. The set
of minima is invariant to permutations and to changes of the signs. Also, according to Proposition~\ref{R-prop}, $\{s^m_i\}$, i.e., the
coordinates of the $\mathbf{s}^m$ components of the \RISA task belong to the set of the minima. \qed
\end{proof}

\subsection{The \CISA Separation Theorem} \label{subsec:CISA-sep-theorem}
The proof of the complex case is similar to the proof of the real case. The difference is in the EPI-type relations
that we apply. Procedure: We define a \C-EPI property and then a \C-w-EPI relation starting from the vector variant of
the \R-EPI relation. Then the proof relies on analogous steps with the real case, that we detail here for the sake of
completeness.

\subsubsection{EPI-type Relations (Complex Case)}\label{subsubsec:CEPIineq-s} Let us consider the vector variant of the \R-EPI relation.
\begin{lemma}[vector-EPI]
For independent (finite covariance random variables) $\mathbf{u}_1,\ldots,\mathbf{u}_L\in\R^q$ holds
\cite{taneja01generalized} that
\begin{equation}
     e^{2H\left(\sum_{i=1}^L\mathbf{u}_i\right)/q}\ge \sum_{i=1}^Le^{2H(\mathbf{u}_i)/q}.\label{eq:v-EPI}
\end{equation}
\end{lemma}

Let us define a similar property for complex random variables:
\begin{definition}[$\C$-EPI]
We say that random variables $u_1,\ldots,u_L\in\C$ satisfy relation $\C$-EPI if
\begin{equation}
    e^{H\left(\sum_{i=1}^Lu_i\right)}\ge \sum_{i=1}^L e ^{H\left(u_i\right)}.\label{eq:C-EPI}
\end{equation}
\end{definition}

\begin{note}
This holds for independent random variables $u_1,\ldots,u_L\in\C$, because according to vector-EPI ($q=2$)
\begin{equation}
    e^{2H\left(\sum_{i=1}^Lu_i\right)/2}\ge \sum_{i=1}^L e ^{2H\left(u_i\right)/2}.
\end{equation}
\end{note}

We need to following lemma:

\begin{lemma}\label{lem:C-suff}
Let us assume that random variables $u_1,\ldots,u_L\in\C$ satisfy condition
\begin{equation}
e^{H(\sum_{i=1}^Lw_iu_i)}\ge\sum_{i=1}^Le^{H(w_iu_i)}\quad \forall \mathbf{w}=[w_1;\ldots;w_L]\in
S^L(\C)\label{eq:C-w-EPI}
\end{equation}
that we shall call condition $\C$-w-EPI. Here, $S^L(\C)$ denotes the \mbox{$L$-dimensional} complex unit sphere, that
is
\begin{equation}
S^L(\C):=\left\{\mathbf{w}=[w_1;\ldots;w_L]\in\C^L:\sum_{i=1}^L\left|w_i\right|^2=1\right\} .
\end{equation}
Then
\begin{equation}
    H\left(\sum_{i=1}^Lw_iu_i\right)\ge\sum_{i=1}^L|w_i|^2H(u_i)\quad \forall \mathbf{w}\in S^L(\C).\label{eq:C-suff}
\end{equation}
\end{lemma}

\begin{proof}
Assume that $\mathbf{w}\in S^L(\C)$. Applying $\ln$ on condition \eqref{eq:C-w-EPI}, and using the monotonicity of the
$\ln$ function, we can see that the first inequality is valid in the following inequality chain
\begin{equation}
H\left(\sum_{i=1}^Lw_iu_i\right)\ge \ln\left(\sum_{i=1}^L e^{H(w_iu_i)}\right)=\ln\left(\sum_{i=1}^Le^{H(u_i)}\cdot
|w_i|^2\right)\ge\sum_{i=1}^L|w_i|^2\cdot\ln\left(e^{H(u_i)}\right)=\sum_{i=1}^L|w_i|^2\cdot H(u_i).
\end{equation}
Then,
\begin{enumerate}
    \item
        we used the relation:
        \begin{equation}
            H(wu)=H(u)+\ln\left(\left|w\right|^2\right) \quad(w,u\in\C)
        \end{equation}
        for the entropy of the transformed variable (see Lemma~\ref{lem:C-H-trafo}). Hence
        \begin{equation}
            e^{H(w_iu_i)}=e^{H(u_i)+\ln\left(\left|w_i\right|^2\right)}=e^{H(u_i)}\cdot
            e^{\ln\left(\left|w_i\right|^2\right)}=e^{H(u_i)}\cdot
            |w_i|^2.\label{eq:C-2entr-transf}
        \end{equation}
    \item
        In the second inequality, we exploited the concavity of $\ln$.\qed
\end{enumerate}
\end{proof}

\subsubsection{Connection to the Cost Function of the \CICA Task}\label{subsubsec:connection2CICA} Now we shall use Lemma~\ref{lem:C-suff} to proceed. The
\CISA Separation Theorem will be a corollary of the following claim:
\begin{proposition}\label{C-prop}
Let $\mathbf{y}=\left[\mathbf{y}^1;\ldots;\mathbf{y}^M\right]=\mathbf{y}(\mathbf{W})=\mathbf{W}\mathbf{s}$, where
$\mathbf{W}\in \mathcal{U}^D$, $\mathbf{y}^m$ is the estimation of the $m^{th}$ component of the \CISA task. Let
$y^m_i$ be the $i^{th}$ complex coordinate of the $m^{th}$ component. Similarly, let $s^m_i$ stand for the $i^{th}$
coordinate of the $m^{th}$ source. Let us assume that the $\mathbf{s}^m$ sources satisfy condition~\eqref{eq:C-suff}.
Then
\begin{equation}\label{eq:C-main-prop}
\sum_{m=1}^M\sum_{i=1}^dH\left(y^m_i\right)\ge \sum_{m=1}^M\sum_{i=1}^dH\left(s^m_i\right).
\end{equation}
\end{proposition}

\begin{proof}
Let us denote the $(i,j)^{th}$ element of matrix $\mathbf{W}$ by $W_{i,j}$. Coordinates of $\mathbf{y}$ and
$\mathbf{s}$ will be denoted by $y_i$ and $s_i$, respectively. Let $\G^1, \ldots, \G^M$ denote the indices belonging to
the $1^{st}, \ldots, M^{th}$. subspaces, that is, $\G^1:=\{1,\ldots,d\},\ldots,\G^M:=\{D-d+1,\ldots,D\}$. Now, writing
the elements of the $i^{th}$ row of matrix multiplication $\mathbf{y}=\mathbf{W}\mathbf{s}$, we have
\begin{equation}
y_i=\sum_{j\in \G^1} W_{i,j}s_j+\ldots+\sum_{j\in \G^M} W_{i,j}s_j\label{eq:C-y=Ws}
\end{equation}
and thus,
\begin{eqnarray}
\lefteqn{H\left(y_i\right)=}\nonumber\\
&&=H\left(\sum_{j\in \G^1} W_{i,j}s_j+\ldots+\sum_{j\in \G^M} W_{i,j}s_j\right)\label{eq:C-H(y=Ws)}\\
&&=H\left(\left(\sum_{l\in\G^1}|W_{i,l}|^2\right)^{\frac{1}{2}}\frac{\sum_{j\in\G^1}W_{i,j}s_j}{\left(\sum_{l\in\G^1}|W_{i,l}|^2\right)^{\frac{1}{2}}}
+ \ldots + \left(\sum_{l\in\G^M}|W_{i,l}|^2\right)^{\frac{1}{2}}\frac{\sum_{j\in\G^M}W_{i,j}s_j}{\left(\sum_{l\in\G^M}|W_{i,l}|^2\right)^{\frac{1}{2}}}\right)\label{eq:C-w^2-in}\\
&&\ge\left(\sum_{l\in\G^1}|W_{i,l}|^2\right)H\left(\frac{\sum_{j\in\G^1}W_{i,j}s_j}{\left(\sum_{l\in\G^1}|W_{i,l}|^2\right)^{\frac{1}{2}}}\right)
+ \ldots    + \left(\sum_{l\in\G^M}|W_{i,l}|^2\right) H\left(\frac{\sum_{j\in\G^M}W_{i,j}s_j}{\left(\sum_{l\in\G^M}|W_{i,l}|^2\right)^{\frac{1}{2}}}\right)\label{eq:C-Lem2-applied}\\
&&=\left(\sum_{l\in\G^1}|W_{i,l}|^2\right)
H\left(\sum_{j\in\G^1}\frac{W_{i,j}}{\left(\sum_{l\in\G^1}|W_{i,l}|^2\right)^{\frac{1}{2}}}s_j\right)
+ \ldots + \left(\sum_{l\in\G^M}|W_{i,l}|^2\right) H\left(\sum_{j\in\G^M}\frac{W_{i,j}}{\left(\sum_{l\in\G^M}|W_{i,l}|^2\right)^{\frac{1}{2}}}s_j\right)\label{eq:C-Lem2-applied-again-pre}\\
&&\ge\left(\sum_{l\in\G^1}|W_{i,l}|^2\right)
\sum_{j\in\G^1}\left|\frac{W_{i,j}}{\left(\sum_{l\in\G^1}|W_{i,l}|^2\right)^{\frac{1}{2}}}\right|^2H\left(s_j\right)\label{eq:C-Lem2-applied-again}
+ \ldots + \left(\sum_{l\in\G^M}|W_{i,l}|^2\right) \sum_{j\in\G^M}\left|\frac{W_{i,j}}{\left(\sum_{l\in\G^M}|W_{i,l}|^2\right)^{\frac{1}{2}}}\right|^2H\left(s_j\right)\\
&&=\sum_{j\in\G^1}|W_{i,j}|^2H\left(s_j\right)+\ldots+\sum_{j\in\G^M}|W_{i,j}|^2H\left(s_j\right)\label{eq:C-H(yi-last)}
\end{eqnarray}
The above steps can be justified as follows:
\begin{enumerate}
    \item
        \eqref{eq:C-H(y=Ws)}: Eq.~\eqref{eq:C-y=Ws} was inserted into the argument of $H$.
    \item
        \eqref{eq:C-w^2-in}: New terms were added for Lemma~\ref{lem:C-suff}.
    \item
        \eqref{eq:C-Lem2-applied}: Sources $\mathbf{s}^m$ are independent of each other and this independence is preserved upon
        mixing \emph{within} the subspaces, and we could also use Lemma~\ref{lem:C-suff}, because $\mathbf{W}$
        is a unitary matrix.
    \item
        \eqref{eq:C-Lem2-applied-again-pre}: Nominators were transferred into the $\sum_j$ terms.
    \item
        \eqref{eq:C-Lem2-applied-again}: Variables $\mathbf{s}^m$ satisfy condition~\eqref{eq:C-suff} according
        to our assumptions.
    \item
        \eqref{eq:C-H(yi-last)}: We simplified the expression after squaring.
\end{enumerate}
Using this inequality, summing it for $i$, exchanging the order of the sums, and making use of the unitary property of
matrix $\mathbf{W}$, we have
\begin{eqnarray}
\sum_{i=1}^DH(y_i)&\ge&\sum_{i=1}^D\left(\sum_{j\in\G^1}|W_{i,j}|^2H\left(s_j\right)+\ldots+\sum_{j\in\G^M}|W_{i,j}|^2H\left(s_j\right)\right)\\
&=&\sum_{j\in\G^1}\left(\sum_{i=1}^D|W_{i,j}|^2\right)H\left(s_j\right)+\ldots+\sum_{j\in\G^M}\left(\sum_{i=1}^D|W_{i,j}|^2\right)H\left(s_j\right)\\
&=&\sum_{j=1}^DH(s_j).
\end{eqnarray}
\qed
\end{proof}

\begin{note}
The proof of the proposition is similar when the dimensions of the subspaces are not constrained to be equal. The
situation is the same in the next theorem.
\end{note}

\subsubsection{Proof of the $\C$-ISA Separation Theorem}\label{subsec:CISA-sepT-proof}

Having this proposition, now we present our main theorem.
\begin{theorem}[Separation Theorem for $\C$-ISA]
Presume that the $\mathbf{s}^m$ sources of the $\C$-ISA model satisfy condition~\eqref{eq:C-suff}, and that
$J(\mathbf{W})=\mbox{$\sum_{m=1}^M\sum_{i=1}^dH\left(y^m_i\right)$}, (\mathbf{W}\in\mathscr{U}^D)$, i.e., the \CICA
cost function has minimum. Then it is sufficient to search for the minimum of the \CISA task
($\mathbf{W}_{\mathrm{\text{\CISA}}}$) as a permutation of the solution of the \CICA task
($\mathbf{W}_{\mathrm{\text{\CICA}}}$). That it, it is sufficient to search in the form
\begin{equation}
    \mathbf{W}_{\mathrm{\text{\CISA}}}=\mathbf{P}\mathbf{W}_{\mathrm{\text{\CICA}}},
\end{equation}
where $\mathbf{P} \left(\in\R^{D\times D}\right)$ is the permutation matrix to be determined.
\end{theorem}
\begin{proof}
\CICA minimizes the l.h.s. of Eq.~\eqref{eq:C-main-prop}, that is, it minimizes \mbox{$\sum_{m=1}^M\sum_{i=1}^dH\left(y^m_i\right)$}. The
set of minima is invariant for permutations and for multiplication of the coordinates by numbers with unit absolute value, and according to
Proposition~\ref{C-prop} $\{s^m_i\}$ (i.e., the coordinates of the $\C$-ISA task) is among the minima.

We can disregard multiplications with unit absolute values, because unitary ambiguity within subspaces are present in
the \CISA task. \qed
\end{proof}

\section{Sufficient Conditions of the Separation Theorem}\label{sec:suff-cond}
In the Separation Theorem, we assumed that relations~\eqref{eq:suff} and \eqref{eq:C-suff} are fulfilled for the
$\mathbf{s}^m$ sources in the real and complex cases, respectively. Here, we shall provide sufficient conditions when
these inequalities are fulfilled.

\subsection{Real Case}

\subsubsection{\R-w-EPI} According to Lemma~\ref{lem:suff}, if the \R-w-EPI property [i.e., \eqref{eq:w-EPI}] holds for sources $\mathbf{s}^m$, then
inequality \eqref{eq:suff} holds, too.

\subsubsection{Real Spherically Symmetric Sources}
\begin{definition}[real spherically symmetric variable]
A random variable \mbox{$\mathbf{u}\in\R^d$} is called real spherically symmetric (or shortly \R-spherical), if
        its density function is not modified by any rotation. Formally, if
        \begin{equation}
        \mathbf{u} \stackrel{\mathrm{distr}}{=}
        \mathbf{O}\mathbf{u},\quad\forall \mathbf{O}\in \mathcal{O}^d,\label{eq:ss-orthinv-def}
        \end{equation}
        where $\stackrel{\mathrm{distr}}{=}$ denotes equality in distribution.
\end{definition}
A \R-spherical random variable has a density function (under mild conditions) and this density function takes constant values on concentric
spheres around the origin. We shall make use of the following well-known properties of spherically symmetric variables
\cite{fang90symmetric,frahm04generalized}:

\begin{lemma}[Identical distribution of 1-dimensional projections - Real case]\label{lem:sph-sym-proj}
Let $\mathbf{v}$ denote a $d$-dimensional variable, which is \R-spherically symmetric. Then the projection of
$\mathbf{v}$ onto lines through the origin have identical univariate distribution.
\end{lemma}

\begin{lemma}[Momenta - Real case]\label{lem:R-SS-momenta}
The expectation value and the variance of a $d$-dimensional $\mathbf{v}$ \R-spherically symmetric variable are
\begin{eqnarray}
            E[\mathbf{v}]&=&\mathbf{0},\label{eq:SS:E}\\
            cov[\mathbf{v}]&=&c(onstant)\cdot \mathbf{I}_d.\label{eq:SS:Var}
\end{eqnarray}
\end{lemma}

Now we are ready to claim the following theorem.
\begin{proposition}
For spherically symmetric sources  $\mathbf{s}^m$ ($m=1,\ldots,M$) with finite covariance Eq.~\eqref{eq:suff} holds.
Further, the stronger \mbox{\R-w-EPI} property [Eq.~\eqref{eq:w-EPI}] also holds and with equality between the two
sides [$\forall \mathbf{w}\in S^d(\R)$].
\end{proposition}
\begin{proof}
Here, we show that the \R-w-EPI property is fulfilled with equality for \R-spherical sources.  According to
\eqref{eq:SS:E}--\eqref{eq:SS:Var}, spherically symmetric sources $\mathbf{s}^m$ have zero expectation values and up to
a constant multiplier they also have identity covariance matrices:
        \begin{eqnarray}
            E[\mathbf{s}^m]&=&\mathbf{0},\\
            cov[\mathbf{s}^m]&=&c^m\cdot \mathbf{I}_d.
        \end{eqnarray}
Note that our constraint on the \RISA task, namely that covariance matrices of the $\mathbf{s}^m$ sources should be
equal to $\mathbf{I}_d$, is fulfilled up to constant multipliers.

Let $P_{\mathbf{w}}$ denote the projection to straight line with direction $\mathbf{w}\in S^d(\R)$, which crosses the
origin, i.e.,
\begin{equation}
P_{\mathbf{w}}:\R^d \ni\mathbf{u}\mapsto\sum_{i=1}^d w_iu_i\in\R.
\end{equation}

In particular, if $\mathbf{w}$ is chosen as the canonical basis vector $\mathbf{e}_i$ (all components are 0, except the
$i^{th}$ component, which is equal to 1), then
\begin{equation}
P_{\mathbf{e_i}}(\mathbf{u})=u_i.
\end{equation}
In this interpretation \R-w-EPI [see Eq.~\eqref{eq:w-EPI}] is concerned with the entropies of the projections of the
different sources onto straight lines crossing the origin. The l.h.s. projects to $\mathbf{w}$, whereas the r.h.s.
projects to the canonical basis vectors. Let $\mathbf{u}$ denote an arbitrary source, i.e., $\mathbf{u}:=\mathbf{s}^m$.
According to Lemma~\ref{lem:sph-sym-proj}, distribution of the spherical $\mathbf{u}$ is the same for all such
projections and thus its entropy is identical. That is,
\begin{eqnarray}
&&\sum_{i=1}^d w_iu_i \stackrel{\mathrm{distr}}{=} u_1
\stackrel{\mathrm{distr}}{=}\ldots
\stackrel{\mathrm{distr}}{=} u_d,\quad \forall \mathbf{w}\in S^d(\R),\\
&&H\left(\sum_{i=1}^d w_iu_i\right) = H\left(u_1\right) =\ldots = H\left(u_d\right),\quad \forall \mathbf{w}\in
S^d(\R).\label{eq:H-inv}
\end{eqnarray}
Thus:
\begin{itemize}
    \item
        l.h.s. of \R-w-EPI: $e^{2H(u_1)}$.
    \item
        r.h.s. of \R-w-EPI:
        \begin{equation}
            \sum_{i=1}^d e^{2H(w_iu_i)}=\sum_{i=1}^de^{2H(u_i)}\cdot
            w_i^2=e^{2H(u_1)}\sum_{i=1}^dw_i^2=e^{2H(u_1)}\cdot 1=e^{2H(u_1)}
        \end{equation}
        At the first step, we used identity \eqref{eq:2entr-transf} for each of the terms.
        At the second step, \eqref{eq:H-inv} was exploited. Then term $e^{H(u_1)}$ was pulled out and
        we took into account that $\mathbf{w}\in S^d(\R)$.
\end{itemize}
\qed
\end{proof}

\begin{note} We note that sources of spherically symmetric distribution have already been used in the context of \RISA
in \cite{hyvarinen00emergence}. In that work, a generative model was assumed. According to the assumption, the distribution of the norms of
sample projections to the subspaces were independent. This way, the task was restricted to spherically symmetric source distributions,
which is a special case of the general \RISA task.
\end{note}

\begin{note}
Spherical variables as well as their non-degenerate affine transforms, the so called elliptical variables (which are equivalent to
spherical ones from the point of view of \RISA) are thoroughly treated in \cite{fang90symmetric,frahm04generalized}.
\end{note}

\subsubsection{Sources Invariant to $90^{\circ}$ Rotation} In the previous
section, we have seen that random variables with density functions invariant to orthogonal transformations (\R-spherical variables) satisfy
the conditions of the \R-ISA Separation Theorem. For mixtures of \mbox{2-dimensional} components ($d=2$), invariance to $90^{\circ}$
rotation suffices. First, we observe that:

\begin{note}\label{note:ort-id-suff} In the \RISA Separation Theorem, it is sufficient if some orthogonal
transformation of the $\mathbf{s}^m$ sources, $\mathbf{C}^m\mathbf{s}^m$ \mbox{($\mathbf{C}^m\in\mathcal{O}^d$)}
satisfy the
 condition \eqref{eq:suff}. In this case, the $\mathbf{C}^m\mathbf{s}^m$ variables are extracted by the permutation
search after the \RICA transformation. Because the \RISA identification has ambiguities up to orthogonal transformation
in the respective subspaces, this is suitable. In other words, for the \RISA identification the existence of an
Orthonormal Basis (ONB) for each $\mathbf{u}:=\mathbf{s}^m\in\R^d$ components is sufficient, on which the
\begin{equation}
h:\R^d\ni\mathbf{w}\mapsto H[\left<\mathbf{w},\mathbf{u}\right>]
\end{equation}
function takes its minimum. [Here, the $\left<\mathbf{w},\mathbf{u}\right>:=\sum_{i=1}^dw_iu_i$ random variable is the
projection of $\mathbf{u}$ to the direction $\mathbf{w}\in S^d(\R)$.] In this case, the  entropy inequality
\eqref{eq:suff} is met with equality on the elements of the ONB.
\end{note}
Now we present our result concerning to the $d=2$ case.
\begin{proposition}\label{prop:90-rot-inv}
Let us suppose, that the density function $f$ of random variable $\mathbf{u}=(u_1,u_2)(=\mathbf{s}^m)\in\R^2$ exhibits
the invariance
\begin{equation}
f(u_1,u_2)=f(-u_2,u_1)=f(-u_1,-u_2)=f(u_2,-u_1)\quad\left(\forall \mathbf{u}\in\R^2\right),\label{sep:suff}
\end{equation}
that is, it is invariant to $90^\circ$ rotation. If function $h(\mathbf{w})=H[\left<\mathbf{w},\mathbf{u}\right>]$ has
minimum on the set $\{\mathbf{w}\ge\mathbf{0}\}\cap S^2(\R)$, it also has minimum on an ONB. \footnote{Relation
$\mathbf{w}\ge \mathbf{0}$ concerns each coordinates.} Consequently, the \RISA task can be identified by the use of the
\RISA Separation Theorem.
\end{proposition}
\begin{proof}
Let
\begin{equation}
\mathbf{R}:=\left[\begin{array}{cc}0&-1\\1&0\end{array}\right]
\end{equation}
denote the matrix of $90^{\circ}$ ccw rotation. Let $\mathbf{w}\in S^2(\R)$. $\left<\mathbf{w},\mathbf{u}\right>\in\R$
is the projection of variable $\mathbf{u}$ onto $\mathbf{w}$. The value of the density function of the random variable
$\left<\mathbf{w},\mathbf{u}\right>$ in $t\in\R$ (we move $t$ in direction $\mathbf{w}$) can be calculated by
integration starting from the point $\mathbf{w}t$, in direction perpendicular to $\mathbf{w}$
\begin{equation}
        f_{y=y(\mathbf{w})=\left<\mathbf{w},\mathbf{u}\right>}(t)=\int_{\mathbf{w}^\perp}
        f (\mathbf{w}t+\mathbf{z})\mathrm{d}\mathbf{z}.\label{eq:proj-int}
\end{equation}
Using the supposed invariance of $f$ and the relation \eqref{eq:proj-int} we have
\begin{equation}
f_{y(\mathbf{w})}=f_{y(\mathbf{Rw})}=f_{y(\mathbf{R}^2\mathbf{w})}=f_{y(\mathbf{R}^3\mathbf{w})},\label{eq:f-invar}
\end{equation}
where `$=$' denotes the equality of functions. Consequently, it is enough to optimize $h$ on the set $\{\mathbf{w}\ge
\mathbf{0}\}$. Let $\mathbf{w}_{min}$ be the minimum of function $h$ on the set $S^2(\R)\cap\{\mathbf{w}\ge
\mathbf{0}\}$. According to Eq.~\eqref{eq:f-invar}, $h$ takes constant and minimal values in the
\[
\{\mathbf{w}_{min},\mathbf{R}\mathbf{w}_{min},\mathbf{R}^2\mathbf{w}_{min},\mathbf{R}^3\mathbf{w}_{min}\}
\]
points. $\{\mathbf{v}_{min},\mathbf{Rv}_{min}\}$ is a suitable ONB in Note~\ref{note:ort-id-suff}.\qed
\end{proof}

\begin{note}
A special case of the requirement \eqref{sep:suff} is invariance to permutation and sign changes, that is
\begin{equation}
     f(\pm u_1,\pm u_2)=f(\pm u_2,\pm u_1).
\end{equation}
In other words, there exists a function $g:\R^2\rightarrow\R$, which is symmetric in its variables and
\begin{equation}
f(\mathbf{u})=g(|u_1|,|u_2|).
\end{equation}
The domain of Proposition~\eqref{prop:90-rot-inv} includes
\begin{enumerate}
    \item
        the formerly presented \R-spherical variables,
    \item
        or more generally, variables with density function of the form
        \begin{equation}
            f(\mathbf{u})=g\left(\sum_i|u_i|^p\right)\quad(p>0).
        \end{equation}
        In the literature \emph{essentially} these variables are called \emph{$L^p(\R)$-norm sphericals} (for $p>1$). Here,
        we use the \emph{$L^p(\R)$-norm spherical} denomination in a slightly extended way, for $p>0$.
\end{enumerate}
\end{note}

\subsubsection{Takano's Dependency Criterion} We have seen that the \R-w-EPI property is sufficient for the \RISA Separation Theorem. In
\cite{takano95inequalities}, sufficient condition is provided to satisfy the EPI condition. The condition is based on the dependencies of
the variables and it concerns the 2-dimensional case. The constraint of $d=2$ may be generalized to higher dimensions. We are not aware of
such generalizations.

We note, however, that \R-w-EPI requires that \R-EPI be satisfied on the surface of the unit sphere. Thus it is satisfactory to consider
the intersection of the conditions detailed in \cite{takano95inequalities} on surface of the unit sphere.

\subsubsection{Summary of Sufficient Conditions (Real Case)} Here, we summarize the presented sufficient conditions of the \R-ISA Separation Theorem. We
have proven, that the requirement described by Eq.~\eqref{eq:suff} for the $\mathbf{s}^m$ sources is sufficient for the
theorem. This holds if the \eqref{eq:w-EPI} \R-w-EPI condition is fulfilled. The stronger \R-w-EPI is valid for
\begin{enumerate}
    \item
        sources satisfying Takano's weak dependency criterion,
    \item
        \R-spherical sources (with equality),
    \item
        sources invariant to $90^\circ$ rotation (for $d=2$). Specially, (i) variables invariant to permutation and sign changes, and
        (ii)$L^p(\R)$-norm spherical variables belong to this family.
\end{enumerate}
These results are summarized schematically in Table~\ref{tab:suffcond-summary}.
\begin{table}
  \centering
  \caption{Sufficient conditions for the \RISA Separation Theorem.}\label{tab:suffcond-summary}
  \[\xymatrixcolsep{2cm}
  \xymatrix{
  &\txt{invariance to $90^{\circ}$ rotation ($d=2$)}\ar@{=>}[ddd]|-{\text{(with `=' for a suitable ONB)}}\ar[dr]|-{\text{specially}}&\\
  &&\text{invariance to sign and permutation}\ar[d]|-{\text{specially}}\\
  &&\text{$L^p(\R)$-norm spherical ($p>0$)}\\
  \txt{Takano's dependency\\($d=2$)} \ar@{=>}[r] & \text{\R-w-EPI}\ar@{=>}[d] & \txt{\R-spherical symmetry (or elliptical)}\ar@{=>}[l]_-{\text{[with `=' for all $\mathbf{w}\in S^d(\R)$]}} \ar[u]|-{\text{generalization for $d=2$}}\\
  &\txt{Equation~\eqref{eq:suff}: sufficient\\ for the \R-Separation Theorem}& }
  \]
\end{table}

\subsection{Complex Case}
We provide sufficient conditions that fulfill \eqref{eq:C-suff}.

\subsubsection{\C-w-EPI} According to Lemma~\ref{lem:C-suff}, if the \C-w-EPI property [i.e., \eqref{eq:C-w-EPI}] holds for sources
$\mathbf{s}^m$, then inequality \eqref{eq:C-suff} holds, too.

\subsubsection{Complex Spherically Symmetric Sources} A complex random variable is complex spherically symmetric, or \C-spherical, for
short, if its density function -- which exists under mild conditions -- is constant on concentric complex spheres. We shall show that
\eqref{eq:C-suff} as well as the stronger \eqref{eq:C-w-EPI} \C-w-EPI relations are fulfilled. We need certain definitions and some basic
features to prove the above statement. Thus, below we shall elaborate on complex sphericals \cite{krishnaiah86complex}.

\begin{definition}[$\C$-spherical variable]
A random variable $\mathbf{v}\in\C^d$ is called \emph{$\C$-spherical}, if $\mathbf{u}=\pv(\mathbf{v})\in\R^{2d}$
\mbox{\R-spherical} \cite{krishnaiah86complex}. Equivalent definition for \C-sphericals is that they are invariant to
unitary transformations. Formally, if
 \begin{equation}
        \mathbf{v} \stackrel{\mathrm{distr}}{=}
        \mathbf{U}\mathbf{v},\quad\forall \mathbf{U}\in \mathcal{U}^d.\label{eq:C-ss-uniinv-def}
        \end{equation}
\end{definition}

We need two basic properties of \C-sphericals \cite{krishnaiah86complex} to prove the theorem. These are analogous to
Lemma~\ref{lem:sph-sym-proj} and Lemma~\ref{lem:R-SS-momenta}:

\begin{lemma}[Identical distribution of 1-dimensional projections - Complex case]
        Projections of $\C$-spherical variables $\mathbf{u}\in\C^d$ onto any unit vectors in $S^d(\C)$ have identical distributions.
        Formally, for $\forall$ $\mathbf{w}_1, \mathbf{w}_2\in S^d(\C)$
        \begin{equation}
            \mathbf{w}_1^*\mathbf{u} \stackrel{\mathrm{distr}}{=} \mathbf{w}_2^*\mathbf{u}\in \C.\label{C-SS-proj-distr}
        \end{equation}
\end{lemma}

\begin{lemma}[Momenta - Complex case]
        For $\C$-spherical variable $\mathbf{u}\in\C^d$:
        \begin{align}
            E[\mathbf{u}]&=\mathbf{0},\label{C-SS:E}\\
            cov[\mathbf{u}]&=c\cdot \mathbf{I}_{d}.\label{C-SS:cov}
        \end{align}
\end{lemma}

We claim the following:

\begin{proposition}
$\C$-spherical sources $\mathbf{s}^m\in \C^d$ ($m=1,\ldots,M$) with finite covariances satisfy condition
\eqref{eq:C-suff} of the $\C$-ISA Separation Theorem. Further, they satisfy $\C$-w-EPI (with equality).
\end{proposition}

\begin{proof}
According to \eqref{C-SS:E} and \eqref{C-SS:cov} for $\C$-spherical components $\mathbf{s}^m\in\C^d$:
$E[\mathbf{s}^m]=\mathbf{0}$, $cov[\mathbf{s}^m]=c^m\cdot \mathbf{I}_d$. Note that our constraint on the $\C$-ISA task,
namely that covariance matrices of the $\mathbf{s}^m$ sources should be equal to identity, is fulfilled up to constant
multipliers.

        Let $P_{\mathbf{w}}$ denote the projection to straight line with direction $\mathbf{w}\in S^d(\C)$, which crosses the origin, i.e.,
        \begin{equation}
            P_{\mathbf{w}}:\C^d \ni\mathbf{u}\mapsto\mathbf{w}^*\cdot\mathbf{u}=\sum_{i=1}^d \bar{w_i}u_i\in\C.
        \end{equation}
        The left and right hand sides of condition \eqref{eq:C-suff} correspond to projection onto vector $\bar{\mathbf{w}}$, and projections onto
        vectors $\mathbf{e}_i=[0;...;0;1;0;...]$ (1 in the $i^{th}$ position and $0$s otherwise), respectively.
        $\bar{\mathbf{w}}\in S^d(\C)\Leftrightarrow\mathbf{w}\in S^d(\C)$, because conjugation preserves length.
        Given property \eqref{C-SS-proj-distr}, the distribution and thus the entropy of these projections are equal.
        That is (let $\mathbf{u}$ denote an arbitrary source, i.e., $\mathbf{u}:=\mathbf{s}^m$),
        \begin{eqnarray}
            &&\sum_{i=1}^d w_iu_i \stackrel{\mathrm{distr}}{=} u_1 \stackrel{\mathrm{distr}}{=}\ldots
            \stackrel{\mathrm{distr}}{=} u_d,\quad \forall \mathbf{w}\in S^d(\C),\\
            &&H\left(\sum_{i=1}^d w_iu_i\right) = H\left(u_1\right) =\ldots = H\left(u_d\right),\quad \forall \mathbf{w}\in S^d(\C).\label{eq:C-H-inv}
        \end{eqnarray}
        Thus:
        \begin{itemize}
            \item
                l.h.s. of $\C$-w-EPI: $e^{H(u_1)}$.
            \item
                r.h.s. of $\C$-w-EPI:
                \begin{equation}
                    \sum_{i=1}^d e^{H(w_iu_i)}=\sum_{i=1}^de^{H(u_i)}\cdot
                    |w_i|^2=e^{H(u_1)}\sum_{i=1}^d|w_i|^2=e^{H(u_1)}\cdot 1=e^{H(u_1)}
                \end{equation}
                At the first step, we used identity \eqref{eq:C-2entr-transf} for each of the terms.
                At the second step, \eqref{eq:C-H-inv} was exploited. Then term $e^{H(u_1)}$ was pulled out and
                we took into account that $\mathbf{w}\in S^d(\C)$.
        \end{itemize}
        \qed
\end{proof}

\subsubsection{Summary of Sufficient Conditions (Complex Case)} Here, we summarize the presented sufficient conditions of the \C-ISA Separation Theorem. We
have proven, that the requirement described by Eq.~\eqref{eq:C-suff} for the $\mathbf{s}^m$ sources is sufficient for
the theorem. This holds if the \eqref{eq:C-w-EPI} \C-w-EPI condition is fulfilled. The stronger \C-w-EPI is valid for
\mbox{\C-spherically} symmetric variables.

These results are summarized schematically in Table~\ref{tab:C-suffcond-summary}.
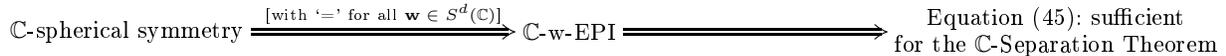
\begin{table}
  \centering
  \caption{Sufficient conditions for the \CISA Separation Theorem.}\label{tab:C-suffcond-summary}
  \[
  \xymatrixcolsep{3.5cm}
  \xymatrix{
  \txt{\C-spherical symmetry}\ar@{=>}[r]^-{\text{[with `=' for all $\mathbf{w}\in S^d(\C)$]}} & \text{\C-w-EPI}\ar@{=>}[r] & \txt{Equation~\eqref{eq:C-suff}: sufficient\\ for the \C-Separation Theorem}}
  \]
\end{table}

\section{Conclusions}\label{sec:conclusion}
In this paper a Separation Theorem,  a decomposition principle, was presented for the \K-Independent Subspace Analysis
(\KISA) problem. If the conditions of the theorem are satisfied then the \KISA task can be solved in 2 steps. The first
step is concerned with the search for 1-dimensional independent components. The second step corresponds to a
combinatorial problem, the search for the optimal permutation. We have shown that spherically symmetric sources (for
the real and the complex cases, too) satisfy the conditions of the theorem. For the real case and for 2-dimensional
sources ($d=2$) invariance to $90^\circ$ rotation, or the Takano's dependency criterion is sufficient for the
separation.

These results underline our experiences that the presented 2 step procedure for solving the \KISA task may produce
higher quality subspaces than sophisticated search algorithms \cite{poczos05independent1}.

Finally we mention that the possibility of this two step procedure (for the real case) was first noted in \cite{cardoso98multidimensional}.

\newpage
\section{Appendix}

\appendix

\section{Uniqueness of \CISA}\label{sec:ISA separability} Here we provide ambiguities of the \CISA task.
The derivation is similar to that of \cite{theis04uniqueness1}, slight modification is used through mappings $\pv,\pM$
[see, Eq.~\eqref{eq:pv} and Eq.~\eqref{eq:pM}].

Notations that we need: Let \mbox{$D=dM$}. Let $Gl(L,\K)$ denote the set of invertible matrices in $\K^{L\times L}$.
Let us decompose matrix $\mathbf{V}\in\C^{D\times D}$ into $d\times d$ blocks:
$\mathbf{V}=\left[\mathbf{V}^{i,j}\right]_{i,j=1,\ldots,M}$ \mbox{$(\mathbf{V}^{i,j}\in\C^{d\times d})$}. We say that
matrix $\mathbf{V}$ is a \emph{$d\times d$ block-permutation matrix}, if there is exactly one index $j$ for $\forall$
$i$ and exactly one $i$ for $\forall$ $j$ ($i,j\in\{1,\ldots,M\}$), that $\mathbf{V}^{ij}\ne\mathbf{0}$, and further,
this block can be inverted. Matrices $\mathbf{B},\mathbf{C}\in\C^{D\times D}$ are \emph{$d$-equivalent} (notation:
$\mathbf{B}\sim_d\mathbf{C}$), if $\mathbf{B}=\mathbf{C}\mathbf{L}$, where $\mathbf{L}\in\C^{D\times D}$ is a $d\times
d$ block-permutation matrix.\footnote{Note: this is an equivalence relation, indeed, because the set of $\mathbf{L}$s
that satisfy the conditions is closed for inversion and for multiplication.} Stochastic variable
\mbox{$\C^D\ni\mathbf{u}=[\mathbf{u}^1;\ldots;\mathbf{u}^M]$} is called \emph{$d$-independent}, if its parts
\mbox{$\mathbf{u}^1,\ldots,\mathbf{u}^M\in\C^d$} are independent. Using $\mathbf{L}$: if $\mathbf{u}$ is
$d$-independent, then $\mathbf{L}\mathbf{u}$ is that, too. Stochastic variable $\mathbf{u}\in\R^L$ is called
\emph{normal}, if every coordinate is normal. Stochastic variable $\mathbf{u}\in\C^L$ is called normal, if both
$\Re(\mathbf{u})$ and $\Im(\mathbf{u})$ are normal. Matrix $\mathbf{B}\in\C^{D\times D}$ is called
\emph{$d$-admissible}, if for
 decomposition $\mathbf{B}=\left[\mathbf{B}^{i,j}\right]_{i,j=1,\ldots,M}$ $(\mathbf{B}^{i,j}\in\C^{d\times d})$ all $\mathbf{B}^{i,j}$ blocks
 are either invertible or indentically 0. (Note: Choosing the coordinates of matrix $\mathbf{B}$ from a continuous distribution, the matrix
is $d$-admissible with probability 1.).

Known properties of $\pM, \pv$  beyond [\eqref{prop:det-trafo}, \eqref{prop:mtx-vec-hom}] \cite{krishnaiah86complex}
are:
\begin{align}
    \pM(\mathbf{M}) \text{ nonsingular (singular)}&\Leftrightarrow \mathbf{M} \text{ nonsingular (singular)},\label{prop:(non)singularity-preserving}\\
    \pv(\mathbf{v}_1+\mathbf{v}_2)&=\pv(\mathbf{v}_1)+\pv(\mathbf{v}_2).\label{prop:vec+-hom}
\end{align}

To prove our statement, we use the following corollary of the Multivariate Skitovitch-Darmois theorem:

\begin{lemma}[Corollary 3.3 in \cite{theis04uniqueness1}]\label{lem:multivariate_S-D}
Let $\mathbf{w}_1=\sum_{m=1}^M\mathbf{B}^m\mathbf{u}^m$ and $\mathbf{w}_2=\sum_{m=1}^M\mathbf{C}^m\mathbf{u}^m$, where
$\mathbf{u}^m$ are independent random variables from $\R^d$, matrices $\mathbf{B}^m,\mathbf{C}^m\in\R^{d\times d}$ are
zeros, or they belong to $GL(d,\R)$. Then, $\mathbf{u}^m$ belonging to $\mathbf{B}^m\mathbf{C}^m\ne\mathbf{0}$ are
normal, provided that $\mathbf{w}_1$ and $\mathbf{w}_2$ are independent.
\end{lemma}

\begin{theorem}[Ambiguities of $\C$-ISA]\label{theorem:C-separability}
Let $Gl(D,\C)\ni\mathbf{B}=[\mathbf{B}^{i,j}]_{i,j=1..M}$ $(\mathbf{B}^{i,j}\in\C^{d\times d})$ $d$-admissible and
\mbox{$\mathbf{s}=[\mathbf{s}^1;\ldots;\mathbf{s}^M]$} $d$-independent $D=dM$-dimensional variable, and none of the
variables $\mathbf{s}^m\in\C^d$ be normal. If $\mathbf{B}\mathbf{s}$ is again $d$-independent, then $\mathbf{B}$ is
d-equivalent to the identity, that is $\mathbf{B}\sim_d\mathbf{I}_D$.
\end{theorem}

\begin{proof}
Indirect. Let us assume that $\mathbf{B}\mathbf{s}$ is $d$-independent, nonetheless $\mathbf{B}\sim_d\mathbf{I}_D$ does
not hold. Then there is a column index $j$ and there are row indices $i_1\ne i_2$ for which
$\mathbf{B}^{i_1,j},\mathbf{B}^{i_2,j}\ne\mathbf{0}$ (and because $\mathbf{B}$ $d$-admissible, thus they are
invertible).\footnote{Reasoning: if for all $j$ there is at most one block (submatrix), which is non-zero, then: (a)
for all $j$ there is exactly one block, which is non-zero and then $\mathbf{B}\sim_d\mathbf{I}_D$, which is a
contradiction, or, (b) there is a $j$ index for which submatrix $\mathbf{B}^{i,j}$ has only zeros, and then the
invertibility of $\mathbf{B}$ is not fulfilled.} Let us take the parts that correspond to indices $i_1,i_2$ off from
$\mathbf{B}\mathbf{s}$:
\begin{align}
\C^{d}\ni\mathbf{y}^{i_1}&=\mathbf{B}^{i_1,j}\mathbf{s}^j + \sum_{m\in\{1,\ldots,M\}\backslash j}\mathbf{B}^{i_1,m}\mathbf{s}^m\\
\C^{d}\ni\mathbf{y}^{i_2}&=\mathbf{B}^{i_2,j}\mathbf{s}^j + \sum_{m\in\{1,\ldots,M\}\backslash
j}\mathbf{B}^{i_2,m}\mathbf{s}^m
\end{align}
Applying $\pv$, and using properties \eqref{prop:vec+-hom} and \eqref{prop:mtx-vec-hom} we have:
\begin{align}
\R^{2d}\ni\pv(\mathbf{y}^{i_1})&=\pM(\mathbf{B}^{i_1,j})\pv(\mathbf{s}^j) + \sum_{m\in\{1,\ldots,M\}\backslash j}\pM(\mathbf{B}^{i_1,m})\pv(\mathbf{s}^m)\\
\R^{2d}\ni\pv(\mathbf{y}^{i_2})&=\pM(\mathbf{B}^{i_2,j})\pv(\mathbf{s}^j) + \sum_{m\in\{1,\ldots,M\}\backslash
j}\pM(\mathbf{B}^{i_2,m})\pv(\mathbf{s}^m)
\end{align}
Taking advantage of \eqref{prop:(non)singularity-preserving}: invertibility of $\mathbf{B}^{i_1,j},\mathbf{B}^{i_2,j}$
is inherited to $\pM(\mathbf{B}^{i_1,j}),\pM(\mathbf{B}^{i_2,j})$. Similarly, matrices $\pM(\mathbf{B}^{i,m})$
($i\in\{i_1,i_2\},m\ne j$) are either zero or they are invertible, according to their ancestor $\mathbf{B}^{i,m}$,
whether it is zero or invertible. If $\mathbf{s}^m\in \C^d$ are independent then variables $\pv(\mathbf{s}^m)\in
\R^{2d}$ are also independent. Thus, as a result of Lemma~\ref{lem:multivariate_S-D}, $\pv(\mathbf{s}^j)$ is normal,
meaning -- by definition -- that $\mathbf{s}^j$ is also normal: a contradiction.\qed
\end{proof}

\begin{note}
\cite{eriksson06complex} has shown an interesting result: for the complex case and for $d=1$ (\CICA task) certain
normal sources can be separated. This result (and thus Theorem~\ref{theorem:C-separability}) may be extended to $d>1$,
too, but we are not aware of such generalization.
\end{note}

\bibliographystyle{splncs}

\begin{thebibliography}{10}

\bibitem{jutten91blind}
Jutten, C., Herault, J.:
\newblock Blind separation of sources: An adaptive algorithm based on
  neuromimetic architecture.
\newblock Signal Processing \textbf{24} (1991)  1--10

\bibitem{comon94independent}
Comon, P.:
\newblock Independent component analysis, a new concept?
\newblock Signal Processing \textbf{36} (1994)  287--314

\bibitem{bell95information}
Bell, A.J., Sejnowski, T.J.:
\newblock An information maximisation approach to blind separation and blind
  deconvolution.
\newblock Neural Computation \textbf{7} (1995)  1129--1159

\bibitem{bell97independent}
Bell, A.J., Sejnowski, T.J.:
\newblock {The `independent components' of natural scenes are edge filters}.
\newblock Vision Research \textbf{37} (1997)  3327--3338

\bibitem{hyvarinen99sparse}
Hyv{\"a}rinen, A.:
\newblock Sparse code shrinkage: Denoising of nongaussian data by maximum
  likelihood estimation.
\newblock Neural Computation \textbf{11} (1999)  1739--1768

\bibitem{kiviluoto98independent}
Kiviluoto, K., Oja, E.:
\newblock Independent component analysis for parallel financial time series.
\newblock In: Proceedings of ICONIP'98. Volume~2. (1998)  895--898

\bibitem{makeig96independent}
Makeig, S., Bell, A.J., Jung, T.P., Sejnowski, T.J.:
\newblock Independent component analysis of electroencephalographic data.
\newblock In: Proceedings of NIPS. Volume~8. (1996)  145--151

\bibitem{Vigario98independent}
Vig{\'a}rio, R., Jousmaki, V., Hamalainen, M., Hari, R., Oja, E.:
\newblock Independent component analysis for identification of artifacts in
  magnetoencephalographic recordings.
\newblock In: Proceedings of NIPS. Volume~10. (1997)  229--235

\bibitem{choi05blind}
Choi, S., Cichocki, A., Park, H.M., Lee, S.Y.:
\newblock Blind source separation and independent component analysis.
\newblock Neural Inf. Proc. Letters and Reviews \textbf{6} (2005)  1--57

\bibitem{hyvarinen00emergence}
Hyv{\"a}rinen, A., Hoyer, P.O.:
\newblock Emergence of phase and shift invariant features by decomposition of
  natural images into independent feature subspaces.
\newblock Neural Computation \textbf{12} (2000)  1705--1720

\bibitem{cardoso98multidimensional}
Cardoso, J.:
\newblock Multidimensional independent component analysis.
\newblock In: Proceedings of ICASSP'98, Seattle, WA. (1998)

\bibitem{theis05blind}
Theis, F.J.:
\newblock Blind signal separation into groups of dependent signals using joint
  block diagonalization.
\newblock In: Proc. {ISCAS} 2005, Kobe, Japan (2005)  5878--5881

\bibitem{kim06independent}
Kim, T., Eltoft, T., Lee, T.W.:
\newblock Independent vector analysis: An extension of {ICA} to multivariate
  components.
\newblock In: Independent Component Analysis and Blind Signal Separation.
  Volume 3889 of LNCS., Springer (2006)  165--172

\bibitem{akaho99MICA}
Akaho, S., Kiuchi, Y., Umeyama, S.:
\newblock {MICA}: Multimodal independent component analysis.
\newblock In: Proceedings of IJCNN. (1999)  927--932

\bibitem{vollgraf01multi}
Vollgraf, R., Obermayer, K.:
\newblock Multi-dimensional {ICA} to separate correlated sources.
\newblock In: Proceedings of NIPS. Volume~14. (2001)  993--1000

\bibitem{bach03finding}
Bach, F.R., Jordan, M.I.:
\newblock Finding clusters in independent component analysis.
\newblock In: Proceedings of ICA2003. (2003)  891--896

\bibitem{poczos05independent1}
P{\'o}czos, B., L{\H{o}}rincz, A.:
\newblock Independent subspace analysis using k-nearest neighborhood distances.
\newblock Artificial Neural Networks: Formal Models and their Applications -
  ICANN 2005, pt 2, Proceedings \textbf{3697} (2005)  163--168

\bibitem{poczos05independent2}
P{\'o}czos, B., L{\H{o}}rincz, A.:
\newblock Independent subspace analysis using geodesic spanning trees.
\newblock In: Proc. of Int. Conf. on Machine Learing (ICML). (2005)  673--680

\bibitem{hulle05edgeworth}
\mbox{Van Hulle}, M.M.:
\newblock Edgeworth approximation of multivariate differential entropy.
\newblock Neural Computation \textbf{17} (2005)  1903--1910

\bibitem{anemuller03complex}
Anem{\"u}ller, J., Sejnowski, T.J., Makeig, S.:
\newblock Complex independent component analysis of frequency-domain
  electroencephalographic data.
\newblock Neural Networks \textbf{16} (2003)  1311--1323

\bibitem{back94blind}
Back, A.D., Tsoi, A.C.:
\newblock Blind deconvolution of signals using a complex recurrent network.
\newblock In: Neural Networks for Signal Processing 4, Proc. of the 1994 IEEE
  Workshop, IEEE Press (1994)  565--574

\bibitem{fiori00blind}
Fiori, S.:
\newblock Blind separation of circularly distributed sources by neural extended
  {APEX} algorithm.
\newblock Neurocomputing Letters \textbf{34} (2000)  239--252

\bibitem{bingham00fast}
Bingham, E., Hyv{\"a}rinen, A.:
\newblock A fast fixed-point algorithm for independent component analysis of
  complex-valued signals.
\newblock International Journal of Neural Systems \textbf{10} (2000)  1--8

\bibitem{cardoso93blind}
Cardoso, J.F., Souloumiac, A.:
\newblock Blind beamforming for non-gaussian signals.
\newblock IEE-Proceedings-F \textbf{140} (1993)  362--370

\bibitem{cardoso99highorder}
Cardoso, J.F.:
\newblock High-order contrasts for independent component analysis.
\newblock Neural Computation \textbf{11} (1999)  157--192

\bibitem{moreau01generalization}
Moreau, E.:
\newblock An any order generalization of {JADE} for complex source signals.
\newblock In: Proc. IEEE International Conference on Acoustics, Speech and
  Signal Processing (ICASSP'2001). Volume~5., Salk Lake City, Utah, USA (2001)
  2805--2808

\bibitem{cardoso96equivariant}
Cardoso, J.F., Laheld, B.:
\newblock Equivariant adaptive source separation.
\newblock IEEE Transactions on Signal Processing (1996)

\bibitem{fiori02complex}
Fiori, S.:
\newblock Complex-weighted one-unit 'rigid-bodies' learning rule for
  independent component analysis.
\newblock Neural Processing Letters \textbf{15} (2002)  275--282

\bibitem{belouchrani97blind}
Belouchrani, A., Abed-Meraim, K., Cardoso, J., Moulines, E.:
\newblock A blind source separation technique using second-order statistics.
\newblock {IEEE} Transactions on Signal Processing \textbf{45} (1997)  434--444

\bibitem{arie00blind}
Yeredor, A.:
\newblock Blind separation of gaussian sources via second-order statistics with
  asymptotically optimal weighting.
\newblock {IEEE} Signal Processing Letters \textbf{7} (2000)

\bibitem{smaragdis98blind}
Smaragdis, P.:
\newblock Blind separation of convolved mixtures in the frequency domain.
\newblock Neurocomputing \textbf{22} (1998)  21--34

\bibitem{calhoun02complex}
Calhoun, V.D., Adali, T.:
\newblock Complex infomax: Convergence and approximation of infomax with
  complex nonlinearities.
\newblock In: NNSP. (2002)

\bibitem{xu04minimax}
Xu, J.W., Erdogmus, D., Rao, Y.N., Pr\'{\i}ncipe, J.C.:
\newblock Minimax mutual information approach for {ICA} of complex-valued
  linear mixtures.
\newblock In: ICA. (2004)  311--318

\bibitem{eriksson04complex}
Eriksson, J., Koivunen, V.:
\newblock Complex-valued {ICA} using second order statistics.
\newblock 2004 {IEEE} Workshop on Machine Learning for Signal Processing (2004)

\bibitem{douglas06equivariant}
Douglas, S.C., Eriksson, J., Koivunen, V.:
\newblock Equivariant algorithms for estimating the strong-uncorrelating
  transform in complex independent component analysis.
\newblock In: Independent Component Analysis and Blind Signal Separation.
  (2006)  57--65

\bibitem{douglas06fixed}
Douglas, S.C., Eriksson, J., Koivunen, V.:
\newblock Fixed-point complex {ICA} algorithms for the blind separation of
  sources using their real or imaginary components.
\newblock In: Independent Component Analysis and Blind Signal Separation.
  (2006)  343--351

\bibitem{poczos05independent3}
P{\'o}czos, B., Tak{\'a}cs, B., L{\H{o}}rincz, A.:
\newblock Independent subspace analysis on innovations.
\newblock Machine Learning: ECML 2005, Proceedings \textbf{3720} (2005)
  698--706

\bibitem{szabo06cross}
Szab{\'o}, Z., P{\'o}czos, B., L{\H{o}}rincz, A.:
\newblock Cross-entropy optimization for independent process analysis.
\newblock In: Independent Component Analysis and Blind Signal Separation.
  Volume 3889 of LNCS., Springer (2006)  909--916

\bibitem{szabo06real}
Szab{\'o}, Z., L{\H{o}}rincz, A.:
\newblock Real and complex independent subspace analysis by generalized
  variance.
\newblock In: ICA Research Network Workshop 2006. (2006)

\bibitem{eriksson06complex}
Eriksson, J.:
\newblock Complex random vectors and {ICA} models: Identifiability, uniqueness
  and separability.
\newblock {IEEE} Transactions on Information Theory \textbf{52} (2006)

\bibitem{eriksson04contributions}
Eriksson, J.:
\newblock Contributions to Theory and Algorithms of Independent Component
  Analysis and Signal Separation.
\newblock PhD thesis, Helsinki University of Technology (2004)

\bibitem{theis04uniqueness1}
Theis, F.:
\newblock Uniqueness of complex and multidimensional independent component
  analysis.
\newblock Signal Processing \textbf{84} (2004)  951--956

\bibitem{theis04uniqueness2}
Theis, F.:
\newblock Uniqueness of real and complex linear independent component analysis
  revisited.
\newblock In: Proc. {EUSIPCO} 2004, Vienna, Austria (2004)  1705--1708

\bibitem{cover91elements}
Cover, T., Thomas, J.:
\newblock Elements of information theory.
\newblock John Wiley and Sons, New York, USA (1991)

\bibitem{krishnaiah86complex}
Krishnaiah, P., Lin, J.:
\newblock Complex elliptically symmetric distributions.
\newblock Communications in Statistics \textbf{15} (1986)  3693--3718

\bibitem{taneja01generalized}
Taneja, I.J.:
\newblock Generalized Information Measures and Their Applications.
\newblock on-line book: www.mtm.ufsc.br/~taneja/book/book.html (2001)

\bibitem{fang90symmetric}
Fang, K.T., Kotz, S., Ng, K.W.:
\newblock Symmetric multivariate and related distributions.
\newblock Chapman and Hall (1990)

\bibitem{frahm04generalized}
Frahm, G.:
\newblock Generalized elliptical distributions: Theory and applications.
\newblock PhD thesis, University of Köln (2004)

\bibitem{takano95inequalities}
Takano, S.:
\newblock The inequalities of {F}isher information and entropy power for
  dependent variables.
\newblock Proceedings of the 7th Japan-Russia Symposium on Probability Theory
  and Mathematical Statistics (1995)

\end{thebibliography}

\end{document}